May 15, 2017

# Melzak's formula for arbitrary polynomials

Khristo N. Boyadzhiev
Ohio Northern University
Department of Mathematics
Ada, Ohio 45810, USA
k-boyadzhiev@onu.edu

**Abstract**. We extend the well-known Melzak binomial transform formula to polynomials of any degree and show some applications.

## 1. Introduction

Let $n \geq 0$ be an integer and let $f(x)$ be a polynomial of degree $\leq n$. Melzak's formula is the following binomial transformation

$$(1) \qquad \sum_{k=0}^{n} \binom{n}{k}(-1)^k \frac{f(x-k)}{y+k} = \frac{n! f(x+y)}{y(y+1)\ldots(y+n)}$$

for any $y \neq 0, -1, -2, \ldots, -n$.

This formula took its name from [5], although it appeared much earlier on p.25 in Nielsen's fundamental study [6]. After the publication of [5] Melzak's formula attracted much attention and was studied and discussed by several authors; see Gould [1], [2], [3], [8], and Prodinger [7]. Gould obtained some interesting extensions, including an extension for polynomials of degree $n+1$, that is, $f(x) = a_{n+1}x^{n+1} + a_n x^n + \ldots$ . In this case

(2) $$\sum_{k=0}^{n}\binom{n}{k}(-1)^k \frac{f(x-k)}{y+k} = \frac{n!f(x+y)}{y(y+1)...(y+n)} - n!a_{n+1}.$$

Melzak's formula is very useful for computing binomial transforms and obtaining binomial identities. For example, when $f=1$ in (1) we find

(3) $$\sum_{k=0}^{n}\binom{n}{k}(-1)^k \frac{1}{y+k} = \frac{n!}{y(y+1)(y+2)...(y+n)} = \frac{1}{y}\binom{n+y}{n}^{-1},$$

or,

(4) $$\sum_{k=0}^{n}\binom{n}{k}(-1)^k \frac{y}{y+k} = \binom{n+y}{n}^{-1},$$

and by inversion of the binomial transform,

(5) $$\sum_{k=0}^{n}\binom{n}{k}(-1)^k \binom{k+y}{k}^{-1} = \frac{y}{y+k}.$$

When $y=1$ in (3) we have

(6) $$\sum_{k=0}^{n}\binom{n}{k}(-1)^k \frac{1}{k+1} = \frac{1}{n+1},$$

or, starting the summation from $k=1$

(7) $$\sum_{k=1}^{n}\binom{n}{k}(-1)^{k-1} \frac{1}{k+1} = \frac{n}{n+1}.$$

## 2. Melzak's formula for arbitrary polynomials

The restriction $\deg f \leq n$ breaks the symmetry in the binomial transform, so it is important to drop this restriction. We present here an extension of (1) to polynomials of any degree (in fact, to any formal power series). We give a short and simple proof of this extension and show that the Stirling numbers of the second kind $S(m,n)$ appear naturally in the new formula.

Before proceeding further, we want to rewrite (1) in a more convenient equivalent form. First, without loss of generality we may assume that $x = 0$. Next, using the well-known partial fraction decomposition

$$(8) \qquad \frac{n!}{y(y+1)...(y+n)} = \sum_{k=0}^{n} \binom{n}{k}(-1)^k \frac{1}{y+k}$$

we write (1) in the form

$$\sum_{k=0}^{n} \binom{n}{k}(-1)^k \frac{f(-k)}{y+k} = f(y)\sum_{k=0}^{n}\binom{n}{k}(-1)^k \frac{1}{y+k}$$

or, equivalently, in the form

$$(9) \qquad \sum_{k=0}^{n} \binom{n}{k}(-1)^k \frac{f(k)}{k-\lambda} = f(\lambda)\sum_{k=0}^{n}\binom{n}{k}(-1)^k \frac{1}{k-\lambda}$$

replacing $y$ by $-\lambda$ and $f(t)$ by $f(-t)$.

**Theorem 1.** Let $f(t) = a_0 + a_1 t + ...$ be a polynomial (or formal power series). Then for any positive integer $n$ and for any complex number $\lambda \neq 0, 1, 2, ...$, we have

$$(10) \qquad \sum_{k=0}^{n}\binom{n}{k}(-1)^k \frac{f(k)}{k-\lambda} = f(\lambda)\sum_{k=0}^{n}\binom{n}{k}\frac{(-1)^k}{k-\lambda} + (-1)^n n! \sum_{m=n+1}^{\infty} a_m \left\{ \sum_{j=0}^{m-1} \lambda^j S(m-j-1,n) \right\}$$

(when $f(t)$ is constant, the sum with Stirling numbers on the RHS is missing).

*Proof.* We shall need the Stirling numbers of the second kind, $S(p,n)$ (see [5]). For every two nonnegative integers $p, n$ we have

$$(11) \qquad (-1)^n n! S(p,n) = \sum_{k=0}^{n}\binom{n}{k}(-1)^k k^p \ .$$

It is well-known that $S(n,n)=1$, $S(n+1,n)=\dfrac{n(n+1)}{2}$, and $S(p,n)=0$ when $p<n$.

Let now $m$ be a positive integer. Obviously,

$$\frac{k^m}{k-\lambda}=\frac{k^m-\lambda^m+\lambda^m}{k-\lambda}=k^{m-1}+\lambda k^{m-2}+\ldots+\lambda^{m-1}+\frac{\lambda^m}{k-\lambda}.$$

Applying the binomial transform we obtain

$$\sum_{k=0}^{n}\binom{n}{k}(-1)^k\frac{k^m}{k-\lambda}=\sum_{p=0}^{m-1}\lambda^p\left\{\sum_{k=0}^{n}\binom{n}{k}(-1)^k k^{m-p-1}\right\}+\lambda^m\sum_{k=0}^{n}\binom{n}{k}(-1)^k\frac{1}{k-\lambda}.$$

In view of (11) this becomes

$$\sum_{k=0}^{n}\binom{n}{k}(-1)^k\frac{k^m}{k-\lambda}=(-1)^n n!\sum_{p=0}^{m-1}\lambda^p S(m-p-1,n)+\lambda^m\sum_{k=0}^{n}\binom{n}{k}(-1)^k\frac{1}{k-\lambda}.$$

We now multiply both sides by the coefficient $a_m$ and sum for $m=1,2,\ldots$ Adding to both sides also the sum

$$\sum_{k=0}^{n}\binom{n}{k}(-1)^k\frac{a_0}{k-\lambda},$$

we obtain the desired formula, namely,

$$\sum_{k=0}^{n}\binom{n}{k}(-1)^k\frac{f(k)}{k-\lambda}=f(\lambda)\sum_{k=0}^{n}\binom{n}{k}\frac{(-1)^k}{k-\lambda}+(-1)^n n!\sum_{m=n+1}^{\infty}a_m\left\{\sum_{p=0}^{m-1}\lambda^p S(m-p-1,n)\right\}.$$

The summation in the second sum on the RHS starts from $m=n+1$, because when $m<n+1$ we have $m-1<n$ and $S(m-1-p,n)=0$. The proof is completed.

When $\deg f(t)\leq n$, the sum on the RHS in (10) is zero and we have the original formula of Melzak (1). When the degree of $f(t)$ equals $n+1$, then $a_m=0$ for $m>n+1$ and (10) turns into

$$\sum_{k=0}^{n}\binom{n}{k}(-1)^k\frac{f(k)}{k-\lambda}=f(\lambda)\sum_{k=0}^{n}\binom{n}{k}\frac{(-1)^k}{k-\lambda}+(-1)^n n!a_{n+1}$$

in accordance with (2). When the degree of $f(t)$ is $n+2$ we have

$$\sum_{k=0}^{n}\binom{n}{k}(-1)^k \frac{f(k)}{k-\lambda} = f(\lambda)\sum_{k=0}^{n}\binom{n}{k}\frac{(-1)^k}{k-\lambda} + (-1)^n n!\left(a_{n+1} + a_{n+2}\left(\frac{n(n+1)}{2}+\lambda\right)\right)$$

etc.

**Remark 2.** Formula (10) can be viewed from a different perspective by using the representation

(12) $$\sum_{k=0}^{n}\binom{n}{k}\frac{(-1)^k}{k-\lambda} = \frac{(-1)^{n+1}n!}{\lambda(\lambda-1)...(\lambda-n)} = (-1)^{n+1}n!\sum_{m=n}^{\infty}\frac{S(m,n)}{\lambda^{m+1}},$$

where the infinite series converges for $|\lambda| > n$ (see [4, (7.47)]).

**Corollary 3**. For any polynomial or formal power series $f(t) = a_0 + a_1 t + ...$ and every $n = 1, 2, ...$

we have

(13) $$\sum_{k=1}^{n}\binom{n}{k}(-1)^{k-1}\frac{f(k)}{k} = f'(0) + f(0)H_n + (-1)^{n-1} n! \sum_{m=n+1}^{\infty} a_m S(m-1,n)$$

where $H_n$ are the harmonic numbers.

$$H_n = 1 + \frac{1}{2} + ... + \frac{1}{n}; \quad H_0 = 0 \ .$$

In particular, when $\deg f \leq n$,

(14) $$\sum_{k=1}^{n}\binom{n}{k}(-1)^{k-1}\frac{f(k)}{k} = f'(0) + f(0)H_n \ .$$

*Proof.* We consider first the case when $f(t)$ has degree $\leq n$ and apply formula (9). On both

sides we separate the term with $k = 0$ to write

$$\frac{f(0)}{-\lambda} + \sum_{k=1}^{n}\binom{n}{k}(-1)^k \frac{f(k)}{k-\lambda} = \frac{f(\lambda)}{-\lambda} + f(\lambda)\sum_{k=1}^{n}\binom{n}{k}(-1)^k \frac{1}{k-\lambda},$$

or, by bringing the first term to the right hand side and multiplying both sides by $-1$

$$\sum_{k=1}^{n}\binom{n}{k}(-1)^{k-1}\frac{f(k)}{k-\lambda}=\frac{f(\lambda)-f(0)}{\lambda}+f(\lambda)\sum_{k=1}^{n}\binom{n}{k}(-1)^{k-1}\frac{1}{k-\lambda}.$$

Setting $\lambda \to 0$ we find

$$\sum_{k=1}^{n}\binom{n}{k}(-1)^{k-1}\frac{f(k)}{k}=f'(0)+f(0)\sum_{k=1}^{n}\binom{n}{k}(-1)^{k-1}\frac{1}{k},$$

which, in view of the well-known identity

$$\sum_{k=1}^{n}\binom{n}{k}(-1)^{k-1}\frac{1}{k}=H_n$$

becomes (14). The general case follows from (10). The first two sums in (10) we manipulate as above and then with $\lambda \to 0$ the last sum in (10) becomes

$$(-1)^{n-1}n!\sum_{m=n+1}^{\infty}a_m S(m-1,n)$$

and the proof is completed.

Identity (14) can be found in the works of Henry Gould [3], [8].

**Remark 4**. When we apply (14) to the polynomial $f(t)=t^{p+1}$ we find

(15) $$\sum_{k=1}^{n}\binom{n}{k}(-1)^{k-1}k^p=0,$$

for any $p<n$. In the case $p=n$, we obtain the well-known property

(16) $$\sum_{k=1}^{n}\binom{n}{k}(-1)^{k-1}k^n=(-1)^{n-1}n!.$$

For $p\geq n$ we have from (15)

$$\sum_{k=1}^{n}\binom{n}{k}(-1)^{k-1}k^p=(-1)^{n-1}n!S(p,n),$$

which is (11).

**Remark 5.** From (2) we can derive the following identity (cf. [3])

(17) $$\sum_{k=0}^{n}\binom{n}{k}(-1)^k \frac{f(x-k)}{(y+k)(z+k)} = \frac{n!}{z-y}\left\{\frac{f(x+y)}{y(y+1)...(y+n)} - \frac{f(x+z)}{z(z+1)...(z+n)}\right\}$$

for any polynomial $f(x)$ with $\deg f \leq n+1$ and variables $y \neq z$. All we need is the decomposition

$$\frac{1}{(y+k)(z+k)} = \frac{1}{z-y}\left\{\frac{1}{y+k} - \frac{1}{z+k}\right\}.$$

The constant $n!a_{n+1}$ from (2) cancels out.

**Corollary 6.** For $f(t) = a_0 + a_1 t + ...$, every integer $n > 0$ and every $\lambda \neq 0,1,2,...$, we have

(18) $$\sum_{k=0}^{n}\binom{n}{k}(-1)^k \frac{f(k)}{(k-\lambda)^2} = f'(\lambda)\sum_{k=0}^{n}\binom{n}{k}\frac{(-1)^k}{k-\lambda} + f(\lambda)\sum_{k=0}^{n}\binom{n}{k}\frac{(-1)^k}{(k-\lambda)^2}$$

$$+ (-1)^n n! \sum_{m=n+1}^{\infty} a_m \left\{\sum_{j=1}^{m-1} j\lambda^{j-1} S(m-j-1,n)\right\}.$$

This follows from Theorem 1 after differentiation with respect to $\lambda$.

Further differentiation in (10) provides a formula for $\sum_{k=0}^{n}\binom{n}{k}(-1)^k \frac{f(k)}{(k-\lambda)^3}$ etc.